\documentclass[twoside]{article}

\usepackage{amsmath,amsthm}     
\usepackage{graphicx}     
\usepackage{quiver}
\usepackage{lipsum}
\usepackage{amsfonts} 
\usepackage{tikz}
\usetikzlibrary{arrows,decorations,patterns,positioning,automata,shadows,fit,shapes,calc,decorations.markings,backgrounds,scopes,decorations.text}
\tikzset{>=latex}

\usepackage[left=4.7cm, right=4.7cm]{geometry}

\usepackage{fancyhdr}
\renewcommand{\leftmark}{The Pseudoinverse of $A=CR$ is $A^+=R^+C^+$ (?)}
\renewcommand{\rightmark}{Micha\l{} P. Karpowicz and Gilbert Strang}
\newtheorem{theorem}{Theorem}
\newtheorem{lemma}{Lemma}
\begin{document}

\title{\mbox{The Pseudoinverse of $A=CR$ is $A^+=R^+C^+$ (?)}}

\author{
\begin{tabular}{c@{\hspace{50pt}}c}
Micha\l{} P. Karpowicz&Gilbert Strang\\
michal.karpowicz@nask.pl&gilstrang@gmail.com\\
NASK PIB&MIT\\
\end{tabular}
\date{}\\
}                      

\maketitle

\begin{abstract}
\noindent
This paper gives three formulas for the pseudoinverse of a matrix product $A = CR$. The first is sometimes correct, the second is always correct, and the third is almost never correct.  But that third randomized pseudoinverse $A^+_r$ may be very useful when $A$ is a very large matrix.

\begin{enumerate}
\item $A^+ = R^+C^+$ when $A  =  CR$ and $C$ has independent columns and $R$ has independent rows.
\item $A^+ = (C^+CR)^+(CRR^+)^+$ is always correct.
\item $A^+_r = (P^TCR)^+P^TCRQ(CRQ)^+ = A^+$ only when $\mathrm{rank}(P^TA) = \mathrm{rank}(AQ) = \mathrm{rank}(A)$ with $A = CR$.
\end{enumerate}

\end{abstract}
\vspace{10pt}


\noindent
The statement in the title is not generally true, as the following example shows:
\begin{gather*}
C = \left[\begin{array}{rrr}1&0\\\end{array}\right] \  
R = \left[\begin{array}{rrr}1\\1\\\end{array}\right]\  
(CR)^+ = \left[ \begin{array}{r}1\\\end{array} \right] \ 
R^+C^+ = \left[\begin{array}{rrr}\dfrac{1}{2}&\dfrac{1}{2}\\\end{array}\right]
\left[\begin{array}{rrr}1\\0\\\end{array}\right]
=\left[\dfrac{1}{2}\right].
\end{gather*}
Here, $C$ has a full row rank and $R$ has a full column rank. We need to have it the other way around in Theorem~\ref{thm1}. Then 
Theorem~\ref{thm2} will give a formula for the pseudoinverse of $A = CR$ that applies in every case.

When the $m$ by $r$ matrix $C$ has $r$ independent columns (full column rank~$r$), and the $r$ by $n$ matrix $R$ has $r$ independent rows (full row rank $r$), the pseudoinverse $C^+$ is the left inverse of $C$, and the pseudoinverse $R^+$ is the right inverse of $R$\,:

\vspace{-11pt}
$$
C^+=(C^\mathrm{T} C)^{-1}C^\mathrm{T}\text{ \;has\; } C^+C=I_{\text{\small$r$}}\quad\mathrm{and}\quad R^+=R^\mathrm{T}(RR^\mathrm{T})^{-1}\text{ \;has\; } RR^+=I_{\text{\small$r$}}.
$$

\noindent
In this case, the $m$ by $n$ matrix $A=CR$ will also have rank $r$ and the statement correctly claims that its $n$ by $m$ pseudoinverse is $A^+=R^+C^+$.

\vspace{10pt}
\begin{theorem}\label{thm1}
The pseudoinverse $A^+$ of a product $A=CR$ is the product of the pseudoinverses $R^+C^+$\hspace{-3pt} (as for inverses)\hspace{-2pt} when\hspace{-1pt} $C$ has full column rank and $R$ has full row rank $r$.
\end{theorem}
\vspace{10pt}

The simplest proof verifies the four Penrose identities \cite{4} that determine\linebreak the pseudoinverse $A^+$ of any matrix $A$\,:

\vspace{-11pt}
\begin{equation}
\!\!AA^+A=A\qquad\quad A^+AA^+=A^+\qquad (A^+A)^\mathrm{T}=A^+A\qquad (AA^+)^\mathrm{T}=AA^+\!\!
\label{eq:Penrose}
\end{equation}

\noindent
Our goal is a different proof of $A^+=R^+C^+$, starting from first principles. We begin with the ``four fundamental subspaces'' associated with any $m$ by $n$ matrix $A$ of rank $r$. Those are the column space $\hbox{\textbf{C}}$ and nullspace $\hbox{\textbf{N}}$ of $A$ and $A^\mathrm{T}$.

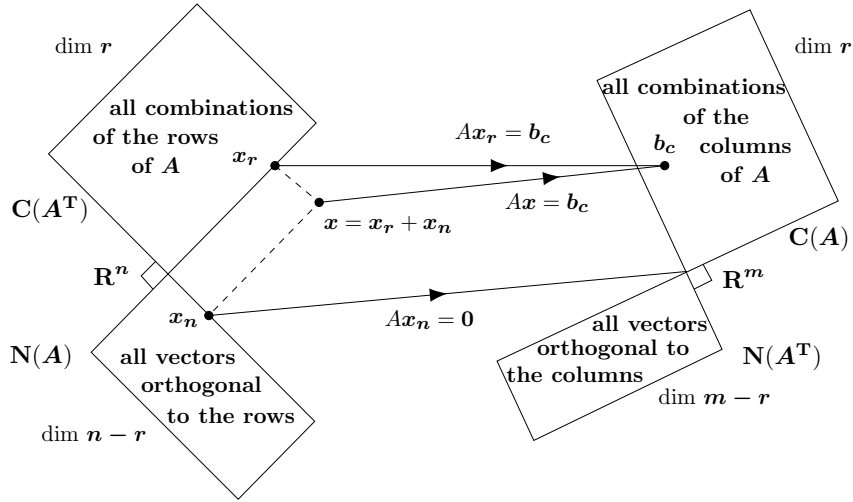
\begin{figure}[!htb]
\centering
 \scalebox{0.85}{
\begin{tikzpicture}[scale=2.3]
\draw (1.07,3.40)--(1.88,2.59)--(0.35,1.03)--(1.35,0.03)--(1.87,0.56)--(0.06,2.38)--cycle;
\draw (0.79,1.65)--(0.69,1.55)--(0.77,1.46);
\node at (0.30,3.12) {dim $\boldsymbol{r}$};
\node at (0.07,2.0) {\scalebox{1.1}{$\hbox{\textbf{C}}(\boldsymbol{A^\mathrm{T}})$}};
\node at (0.37,0.47) {dim $\boldsymbol{n-r}$};
\node at (0.01,1.0) {\scalebox{1.1}{$\hbox{\textbf{N}}(\boldsymbol{A})$}};
\node at (0.49,1.55) {\scalebox{1.1}{$\hbox{\textbf{R}}^{\text{\scriptsize{$\boldsymbol{n}$}}}$}};
\node at (1.1,2.7) {\textbf{all combinations}};
\node at (0.79,2.5) {\textbf{of the rows}};
\node at (0.79,2.3) {\textbf{of $\boldsymbol{A}$}};
\node at (1.4,2.35) {$\boldsymbol{x}_\text{\footnotesize ${\boldsymbol{r}}$}$};
\draw [fill=black](1.6, 2.3) circle (0.8pt);
\draw [fill=black](4.25, 2.3) circle (0.8pt);
\draw [fill=black](1.90,2.05) circle (0.8pt);
\draw [fill=black](1.15,1.28) circle (0.8pt);
\draw [dashed](1.6, 2.3)--(1.90,2.05)--(1.15,1.28);
\draw[decoration={markings, mark=at position 0.6 with {\arrow[scale=2]{>}}},postaction={decorate}](1.6, 2.3)--(4.25, 2.3);
\node at (3.15,2.52) {$A\boldsymbol{x}_\text{\footnotesize ${\boldsymbol{r}}$}=\boldsymbol{b}_\text{\footnotesize ${\boldsymbol{c}}$}$};
\draw[decoration={markings, mark=at position 0.7 with {\arrow[scale=2]{>}}},postaction={decorate}](1.90,2.05)--(4.25, 2.3);
\node at (2.37,1.9) {$\boldsymbol{x}=\boldsymbol{x}_\text{\footnotesize ${\boldsymbol{r}}$}+\boldsymbol{x}_\text{\footnotesize ${\boldsymbol{n}}$}$};
\node at (3.45,2.05) {$A\boldsymbol{x}=\boldsymbol{b}_\text{\footnotesize ${\boldsymbol{c}}$}$};
\draw[decoration={markings, mark=at position 0.5 with {\arrow[scale=2]{>}}},postaction={decorate}](1.15,1.28)--(4.4, 1.58);
\node at (0.98,1.25) {$\boldsymbol{x}_\text{\footnotesize ${\boldsymbol{n}}$}$};
\node at (2.65,1.25) {$A\boldsymbol{x}_\text{\footnotesize ${\boldsymbol{n}}$}=\boldsymbol{0}$};

\node at (0.94,1.0) {\textbf{all vectors}};
\node at (1.10,0.8) {\textbf{orthogonal}};
\node at (1.30,0.6) {\textbf{to the rows}};


\draw (3.79,2.88)--(4.83,3.36)--(5.43,2.06)--(3.11,0.97)--(3.36,0.43)--(4.64,1.05)--cycle;
\draw (4.51,1.63)--(4.57,1.52)--(4.45,1.46);
\node at (5.33,3.12) {dim $\boldsymbol{r}$};
\node at (5.30,1.8) {\scalebox{1.1}{$\hbox{\textbf{C}}(\boldsymbol{A})$}};
\node at (4.59,0.72) {dim $\boldsymbol{m-r}$};
\node at (5.05,1.0) {\scalebox{1.1}{$\hbox{\textbf{N}}(\boldsymbol{A^\mathrm{T}})$}};
\node at (4.79,1.53) {\scalebox{1.1}{$\hbox{\textbf{R}}^{\text{\scriptsize{$\boldsymbol{m}$}}}$}};
\node at (4.45,2.84) {\textbf{all combinations}};
\node at (4.6,2.64) {\textbf{of the}};
\node at (4.8,2.44) {\textbf{columns}};
\node at (4.8,2.24) {\textbf{of} $\boldsymbol{A}$};
\node at (4.25,2.45) {$\boldsymbol{b}_\text{\footnotesize ${\boldsymbol{c}}$}$};%
\node at (4.15,1.23) {\textbf{all vectors}};
\node at (3.9,1.05) {\textbf{orthogonal to}};
\node at (3.64,0.88) {\textbf{the columns}};
\end{tikzpicture}
}
\vspace{-5pt}
\caption{
$A\boldsymbol{x}_\text{\footnotesize ${\boldsymbol{r}}$}=\boldsymbol{b}$ is in the column space of $A$ and $A\boldsymbol{x}_\text{\footnotesize ${\boldsymbol{n}}$}=\boldsymbol{0}$. The complete solution to $A\boldsymbol{x}=\boldsymbol{b}$ is $\boldsymbol{x}=$ \textbf{\,one\;} $\boldsymbol{x}_\text{\footnotesize ${\boldsymbol{r}}$}\,+$ \textbf{\,any\;} $\boldsymbol{x}_\text{\footnotesize ${\boldsymbol{n}}$}$. 
}\label{fig:1}
\end{figure}

Note that every matrix $A$ gives an invertible map (Figure \ref{fig:1}) from its row space $\hbox{\textbf{C}}(A^\mathrm{T})$ to its column space $\hbox{\textbf{C}}(A)$. If $\boldsymbol{x}$ and $\boldsymbol{y}$ are in the row space and $A\boldsymbol{x}=A\boldsymbol{y}$, then $A(\boldsymbol{x}-\boldsymbol{y})=\boldsymbol{0}$. Therefore $\boldsymbol{x}-\boldsymbol{y}$ is in the nullspace $\hbox{\textbf{N}}(A)$ as well as the row space. So $\boldsymbol{x}-\boldsymbol{y}$ is orthogonal to itself and $\boldsymbol{x}=\boldsymbol{y}$.

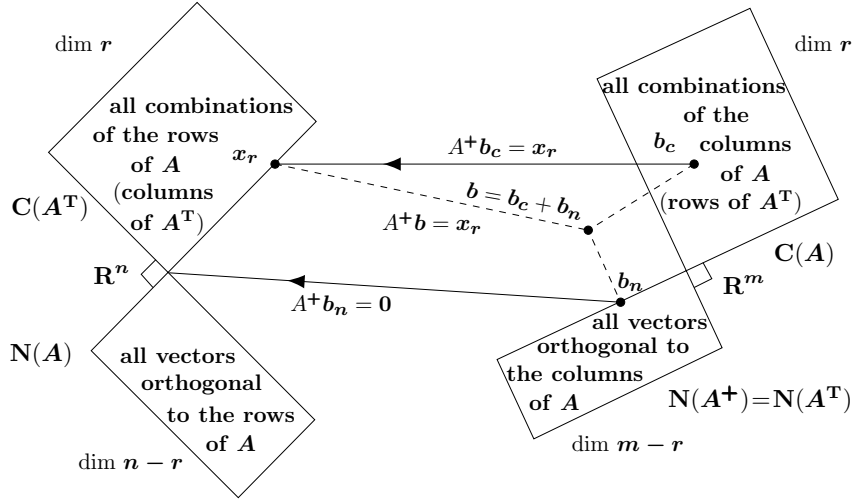
\begin{figure}[th!]
\centering
\scalebox{0.85}{\hspace{9pt}
\begin{tikzpicture}[scale=2.3]
\draw (1.07,3.40)--(1.88,2.59)--(0.35,1.03)--(1.35,0.03)--(1.87,0.56)--(0.06,2.38)--cycle;
\draw (0.79,1.65)--(0.69,1.55)--(0.77,1.46);
\node at (0.30,3.12) {dim $\boldsymbol{r}$};
\node at (0.07,2.0) {\scalebox{1.1}{$\hbox{\textbf{C}}(\boldsymbol{A^\mathrm{T}})$}};
\node at (0.62,0.27) {dim $\boldsymbol{n-r}$};
\node at (0.01,1.0) {\scalebox{1.1}{$\hbox{\textbf{N}}(\boldsymbol{A})$}};
\node at (0.49,1.55) {\scalebox{1.1}{$\hbox{\textbf{R}}^{\text{\scriptsize{$\boldsymbol{n}$}}}$}};
\node at (1.1,2.7) {\textbf{all combinations}};
\node at (0.79,2.5) {\textbf{of the rows}};
\node at (0.79,2.3) {\textbf{of} $\boldsymbol{A}$};
\node at (0.84,2.1) {(\textbf{columns}};
\node at (0.87,1.9) {\textbf{of} $\boldsymbol{A^\mathrm{T}}$)};
\node at (1.4,2.35) {$\boldsymbol{x}_\text{\footnotesize ${\boldsymbol{r}}$}$};
\draw [fill=black](1.6, 2.3) circle (0.8pt);
\draw [fill=black](4.45, 2.3) circle (0.8pt);
\draw [fill=black](3.73,1.85) circle (0.8pt);
\draw [fill=black](3.95,1.36) circle (0.8pt);
\draw [dashed](4.45, 2.3)--(3.73,1.85)--(3.95,1.36);
\draw [dashed](1.6, 2.3)--(3.73,1.85);
\draw[decoration={markings, mark=at position 0.3 with {\arrow[scale=2]{<}}},postaction={decorate}](1.6, 2.3)--(4.45, 2.3);
\node at (3.15,2.42) {$A^{\pmb{+}}\boldsymbol{b}_\text{\footnotesize ${\boldsymbol{c}}$}=\boldsymbol{x}_\text{\footnotesize ${\boldsymbol{r}}$}$};
\node at (3.30,2.05)[rotate=-10] {$\boldsymbol{b}=\boldsymbol{b}_\text{\footnotesize ${\boldsymbol{c}}$}+\boldsymbol{b}_\text{\footnotesize ${\boldsymbol{n}}$}$};
\node at (2.67,1.9) {$A^{\pmb{+}}\boldsymbol{b}=\boldsymbol{x}_\text{\footnotesize ${\boldsymbol{r}}$}$};
\draw[decoration={markings, mark=at position 0.3 with {\arrow[scale=2]{<}}},postaction={decorate}](0.88,1.56)--(3.95,1.36);
\node at (4.02,1.5) {$\boldsymbol{b}_\text{\footnotesize ${\boldsymbol{n}}$}$};
\node at (2.05,1.35) {$A^{\pmb{+}}\boldsymbol{b}_\text{\footnotesize ${\boldsymbol{n}}$}=\boldsymbol{0}$};

\node at (0.94,1.0) {\textbf{all vectors}};
\node at (1.10,0.8) {\textbf{orthogonal}};
\node at (1.30,0.6) {\textbf{to the rows}};
\node at (1.30,0.4) {\textbf{of} $\boldsymbol{A}$};

\draw (3.79,2.88)--(4.83,3.36)--(5.43,2.06)--(3.11,0.97)--(3.36,0.43)--(4.64,1.05)--cycle;
\draw (4.51,1.63)--(4.57,1.52)--(4.45,1.46);
\node at (5.33,3.12) {dim $\boldsymbol{r}$};
\node at (4.0,0.39) {dim $\boldsymbol{m-r}$};
\node at (4.79,1.48) {\scalebox{1.1}{$\hbox{\textbf{R}}^{\text{\scriptsize{$\boldsymbol{m}$}}}$}};
\node at (5.2,1.68) {\scalebox{1.1}{$\hbox{\textbf{C}}(\boldsymbol{A})$}};
\node at (4.9,0.7) {\scalebox{1.1}{$\hbox{\textbf{N}}(\boldsymbol{A^{\pmb{+}}})\hspace{-2pt}=\hspace{-2pt}\hbox{\textbf{N}}(\boldsymbol{A^\mathrm{T}})$}};
\node at (4.45,2.84) {\textbf{all combinations}};
\node at (4.6,2.64) {\textbf{of the}};
\node at (4.8,2.44) {\textbf{columns}};
\node at (4.8,2.24) {\textbf{of} $\boldsymbol{A}$};
\node at (4.7,2.04) {(\textbf{rows of} $\boldsymbol{A^\mathrm{T}}$)};
\node at (4.25,2.45) {$\boldsymbol{b}_\text{\footnotesize ${\boldsymbol{c}}$}$};%
\node at (4.15,1.23) {\textbf{all vectors}};
\node at (3.9,1.05) {\textbf{orthogonal to}};
\node at (3.64,0.88) {\textbf{the columns}};
\node at (3.5,0.68) {\textbf{of} $\boldsymbol{A}$};
\end{tikzpicture}
}
\vspace{-5pt}
\caption{The four subspaces for $A^{\pmb{+}}$ are the four subspaces for $A^\mathrm{T}$.}
  \label{fig:2}
\end{figure}

The pseudoinverse $A^+$ in Figure \ref{fig:2} inverts the row space to column space map in Figure \ref{fig:1} when $\hbox{\textbf{N}}(A^T) = \hbox{\textbf{N}}(A^+)$. If $\boldsymbol{b}_c$ and $\boldsymbol{d}_c$ are in the column space $\hbox{\textbf{C}}(A)$ and $\boldsymbol{x} = A^+\boldsymbol{b}_c=A^+\boldsymbol{d}_c$, then $A^+(\boldsymbol{b}_c-\boldsymbol{d}_c)=\boldsymbol{0}$ and $\boldsymbol{b}_c-\boldsymbol{d}_c$ is in the nullspace $\hbox{\textbf{N}}(A^+)$. Therefore, if $\hbox{\textbf{N}}(A^+) = \hbox{\textbf{N}}(A^T)$, then $\boldsymbol{b}_c-\boldsymbol{d}_c$ is in the nullspace of $A^T$ and the column space of $A$, it is orthogonal to itself and $(\boldsymbol{b}_c-\boldsymbol{d}_c)^T(\boldsymbol{b}_c-\boldsymbol{d}_c) = 0$ with $\boldsymbol{b}_c=\boldsymbol{d}_c$. For all vectors $\boldsymbol{b}_n$ in the orthogonal complement of $\hbox{\textbf{C}}(A)$, we have $A^\mathrm{T}\boldsymbol{b}_n=\boldsymbol{0}$ and  $A^+\boldsymbol{b}_n=\boldsymbol{0}$. Theorem \ref{thm1} shows why $\hbox{\textbf{N}}(A^T) = \hbox{\textbf{N}}(A^+)$ and $A^+$ is the inverse map.

That upper map from the row space $\hbox{\textbf{C}}(A^\mathrm{T})$ to the column space $\hbox{\textbf{C}}(A)$ is inverted by the pseudoinverse $A^+$ in Figure \ref{fig:2}. And the nullspace of $\boldsymbol{A^+}$ is the same as the nullspace of $\boldsymbol{A^\mathrm{T}}$---the orthogonal complement of $\hbox{\textbf{C}}(A)$. \textbf{Thus $\boldsymbol{1/0=0}$ for $\boldsymbol{A^+}$}. In the extreme case, the pseudoinverse of $A=$ zero matrix $(m$ by $n$) is $A^+=$ zero matrix ($n$ by $m$).

\newpage
\noindent
The proof of Theorem 1 begins with a simple Lemma.

\begin{lemma}
We have $\boldsymbol{\mathrm{C}}(A)=\boldsymbol{\mathrm{C}}(AB)$ if and only if $\text{rank}(A)=\text{rank}(AB)$, and $\boldsymbol{\mathrm{N}}(B)=\boldsymbol{\mathrm{N}}(AB)$ if and only if $\text{rank}(B)=\text{rank}(AB)$.
\end{lemma}

\vspace{5pt}
\noindent
\textbf{Proof\,:}\quad Always $\boldsymbol{\mathrm{C}}(A)\supset\boldsymbol{\mathrm{C}}(AB)$ and $\boldsymbol{\mathrm{N}}(AB)\supset\boldsymbol{\mathrm{N}}(B)$. In each pair, equality of dimensions guarantees equality of spaces. The column space can only become smaller, and the nullspace larger.

\vspace{9pt}
\noindent
\textbf{Proof of Theorem 1:}
We will show that $\boldsymbol{\mathrm{N}}(A^\mathrm{T})\subset\boldsymbol{\mathrm{N}}(A^+)$ and $\boldsymbol{\mathrm{N}}(A^\mathrm{T})\supset\boldsymbol{\mathrm{N}}(A^+)$, and that $\boldsymbol{\mathrm{C}}(A^\mathrm{T})\subset\boldsymbol{\mathrm{C}}(A^+)$ and $\boldsymbol{\mathrm{C}}(A^\mathrm{T})\supset\boldsymbol{\mathrm{C}}(A^+)$ when $A^+=R^+C^+$.

\vspace{9pt}
\noindent
\textbf{Step 1\,:}\quad $\hbox{\textbf{N}}(A^\mathrm{T})\subset\hbox{\textbf{N}}(A^+)$. Every vector $\boldsymbol{b_n}$ in the nullspace of $A^\mathrm{T}$ is also in the nullspace of $A^+$.

\vspace{9pt}
\noindent
If $A^\mathrm{T}\boldsymbol{b_n}=0$, then $R^\mathrm{T}C^\mathrm{T}\boldsymbol{b_n}=0$. Since $C^\mathrm{T}$ has full column rank, we conclude that $C^\mathrm{T}\boldsymbol{b_n}=(C^\mathrm{T}C)^{-1}C^\mathrm{T}\boldsymbol{b_n}=C^+\boldsymbol{b_n}=0$ and $\boldsymbol{b_n}$ is in the nullspace of $C^+$. However, by Lemma 1, we also have $\hbox{\textbf{N}}(C^+)=\hbox{\textbf{N}}(R^+C^+)=\hbox{\textbf{N}}(A^+)$. Every $\boldsymbol{b_n}$ in the nullspace of $A^\mathrm{T}$ is also in the nullspace of $A^+$.

\vspace{9pt}
\noindent
\textbf{Step 2\,:}\quad $\hbox{\textbf{N}}(A^\mathrm{T})\supset\hbox{\textbf{N}}(A^+)$. Every vector $\boldsymbol{b_n}$ in the nullspace of $A^+$ is also in the nullspace of $A^\mathrm{T}$.

\vspace{9pt}
\noindent
Suppose that $A^+\boldsymbol{b_n}=R^+C^+\boldsymbol{b_n}=\boldsymbol{0}$. Since $RR^+=I_r$, we have $AA^+\boldsymbol{b_n}=CRR^+C^+\boldsymbol{b_n}=C(C^\mathrm{T}C)^{-1}C^\mathrm{T}\boldsymbol{b_n}=\boldsymbol{0}$. Since $C(C^\mathrm{T}C)^{-1}$ has full column rank, we see that $C^\mathrm{T}\boldsymbol{b_n}=\boldsymbol{0}$, so $\boldsymbol{b_n}$ is orthogonal to every column of $C$. Therefore $\boldsymbol{b_n}\in\hbox{\textbf{N}}(A^\mathrm{T})$.

\vspace{9pt}
\noindent
\textbf{Step 3\,:}\quad $\hbox{\textbf{C}}(A^\mathrm{T})\subset\hbox{\textbf{C}}(A^+)$. Every vector $\boldsymbol{x}$ in the row space of $A$ is also in the column space of $A^+$.

\vspace{9pt}
\noindent
We have $\boldsymbol{x} = A^\mathrm{T}\boldsymbol{y}=R^\mathrm{T}C^\mathrm{T}\boldsymbol{y}=R^\mathrm{T}[(RR^\mathrm{T})^{-1}RR^\mathrm{T}]C^\mathrm{T}\boldsymbol{y}=R^+R(R^\mathrm{T}C^\mathrm{T}\boldsymbol{y})=R^+R\boldsymbol{x}$. So by Lemma 1, we conclude that $\boldsymbol{x}\in\hbox{\textbf{C}}(R^+)=\hbox{\textbf{C}}(R^+C^+)=\hbox{\textbf{C}}(A^+)$.

\vspace{9pt}
\noindent
\textbf{Step 4\,:}\quad $\hbox{\textbf{C}}(A^\mathrm{T})\supset\hbox{\textbf{C}}(A^+)$. Every vector $\boldsymbol{x}$ in the column space of $A^+$ is also in the row space of $A$.

\vspace{9pt}
\noindent
Given $\boldsymbol{b}=A\boldsymbol{x}=CR\boldsymbol{x}$, we have

\vspace{-11pt}
$$\boldsymbol{x}=A^+\boldsymbol{b}=R^+C^+CR\boldsymbol{x}=R^\mathrm{T}(RR^\mathrm{T})^{-1}(C^\mathrm{T}C)^{-1}C^\mathrm{T}CR\boldsymbol{x}=R^\mathrm{T}(RR^\mathrm{T})^{-1}R\boldsymbol{x}$$

\noindent
So again by Lemma 1, we conclude that $\boldsymbol{x}\in\hbox{\textbf{C}}(R^\mathrm{T})=\hbox{\textbf{C}}(R^\mathrm{T}C^\mathrm{T})=\hbox{\textbf{C}}(A^\mathrm{T})$.

\vspace{10pt}
An $m$ by $n$ matrix of rank $r$ has many factorizations. The simplest form of $A=CR$ fills $C$ with the first $r$ independent columns of $A$. Those columns are a~basis for the column space of $A$. Then column $j$ of $R$ specifies the combination of columns of $C$ that produces column $j$ of $A$. This example has column 3 = column 1 + column 2:
$$
A=\left[\begin{array}{rrr}1&4&5\\2&3&5\\\end{array}\right]=\left[\begin{array}{rr}1&4\\2&3\\\end{array}\right]\left[\begin{array}{rrr}1&0&1\\0&1&1\\\end{array}\right]=CR.
$$

\noindent
That matrix $R$ contains the nonzero rows of the \textbf{reduced row echelon form} of $A$. It also reveals nullspace $\hbox{\textbf{N}}(A)=\hbox{\textbf{N}}(R)$. 

$A=CR$ expresses this classical ``elimination by row operations'' as a matrix factorization \cite{5,6}. It is slower and less stable numerically than the SVD, but rational $A$ produces rational $C$ and $R$ and $C^+$ and $R^+$. The following useful formula shows that  explicitly:
\begin{equation}
A^+=R^+C^+=R^\mathrm{T}(RR^\mathrm{T})^{-1}(C^\mathrm{T}C)^{-1}C^\mathrm{T}=R^\mathrm{T}(C^\mathrm{T}AR^\mathrm{T})^{-1}C^\mathrm{T}.
\label{eq:pinverse-formula}
\end{equation}
\noindent
The inverses of square rational matrices $RR^T$ and $C^TC$ are formed by elementary operations and division by a rational determinant. Therefore, the pseudoinverse of $A$ is rational when $A$ is rational.

The ranks of $A^+$ and $A$ are equal when $A^+$ is the inverse map. Then,
$$
\boldsymbol{b} = A\boldsymbol{x}+(I_m-AA^+)\boldsymbol{w}
$$
is a solution of $A^+\boldsymbol{b} = \boldsymbol{x}$ for arbitrary $\boldsymbol{w}$. We have $A^+\boldsymbol{b} = A^+A\boldsymbol{x} = \boldsymbol{x}$, since $A^+AA^+ = A^+$ when $A^+$ and $A$ have equal ranks. 

That is not always the case and not always necessary. Consider the rank-deficient 
$$
A=
\left[\begin{array}{rrr}1&0\\0&0\\0&0\\\end{array}\right]
=
\left[\begin{array}{rrr}1\\0\\0\\\end{array}\right]
\left[\begin{array}{rrr}1&0\\\end{array}\right]
=C_0R_0.
$$
When we complete $C_0$ and $R_0$ to get square invertible $[C_0, C_1]$ and $[R_0,  R_1]^T$, then the generalized $A = CR$ becomes
$$
A= 
\left[\begin{array}{rrr}C_0 & C_1\\\end{array}\right]
\left[\begin{array}{rrr}I_r & 0\\0& 0\\\end{array}\right]
\left[\begin{array}{rrr}R_0 \\R_1\\\end{array}\right] = \bar{C}\bar{R}
$$
and its generalized inverse is
$$
G = \left[\begin{array}{rrr}R_0 \\R_1\\\end{array}\right]^{-1}
\left[\begin{array}{rrr}I_r & Z_{11}\\Z_{21}& Z_{22}\\\end{array}\right]
\left[\begin{array}{rrr}C_0 & C_1\\\end{array}\right]^{-1}
$$
for arbitrary $Z_{ij}$ of proper size \cite{8}. The first Penrose identity holds, $AGA = A$, but the second one gives
$$
GAG = \left[\begin{array}{rrr}R_0 \\R_1\\\end{array}\right]^{-1}
\left[\begin{array}{rrr}I_r & Z_{12}\\Z_{21}& Z_{21}Z_{12}\\\end{array}\right]
\left[\begin{array}{rrr}C_0 & C_1\\\end{array}\right]^{-1}.
$$
We have $GAG=G$ only if $Z_{22} = Z_{21}Z_{12}$. Then $A$ and $G$ have equal ranks and $G$ is the inverse map. Otherwise, the nullspaces $\hbox{\textbf{N}}(A^T)$ and $\hbox{\textbf{N}}(G)$ are different. For
$$
A^T=
\left[\begin{array}{rrr}1&0&0\\0&0&0\\\end{array}\right]
=\left[\begin{array}{rrr}1\\0\\\end{array}\right]
\left[\begin{array}{rrr}1&0&0\\\end{array}\right] = C_2R_2
$$
the nullspace $\hbox{\textbf{N}}(A^T) = \hbox{\textbf{N}}(R_2)$ is spanned by $(0,1,0)^T$ and $(0,0,1)^T$. For:
$$
G=\left[\begin{array}{rrr}1&3&2\\3&3&2\\\end{array}\right]
=
\left[\begin{array}{rrr}1&3\\3&3\\\end{array}\right]
\left[\begin{array}{rrr}1&0&0\\0&1&2/3\\\end{array}\right] = C_3R_3,
$$
we see that $\hbox{\textbf{N}}(G) = \hbox{\textbf{N}}(R_3)$ is spanned by $(0,-2/3,0)^T$. Here $AGA = A$, but $GAG$ is not equal to $G$ with $1 = \mathrm{rank}(A) < \mathrm{rank}(G) = 2$. Also, neither $AG$ nor $GA$ gives an identity matrix. Still, 
$$
\boldsymbol{x} = G\boldsymbol{b}+(I_n-GA)\boldsymbol{z}
$$ 
solves $A\boldsymbol{x} = AG\boldsymbol{b}+(A-AGA)\boldsymbol{z} = AG\boldsymbol{b} =\boldsymbol{b}$ for arbitrary $\boldsymbol{z}$.

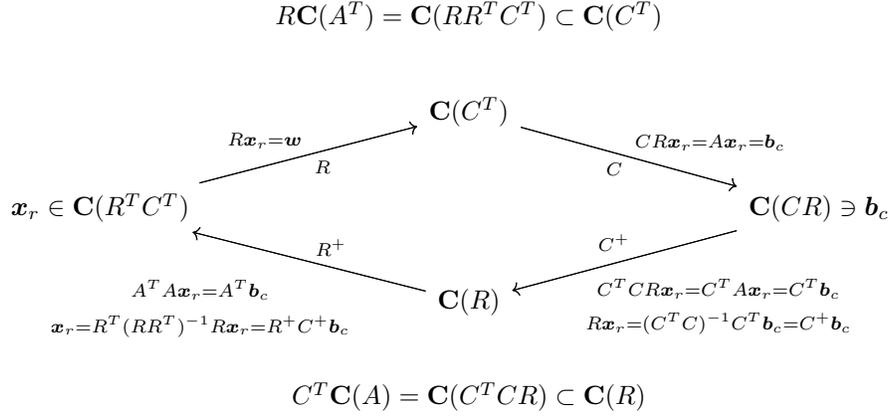
\begin{figure}[!t]
\centering
\[\begin{tikzcd}
& {R\boldsymbol{\mathrm{C}}(A^T) = \boldsymbol{\mathrm{C}}(RR^TC^T) \subset\boldsymbol{\mathrm{C}}(C^T)} \\
& {\hbox{\textbf{C}}(C^T)} \\
{\boldsymbol{x}_r\in\hbox{\textbf{C}}(R^TC^T)} && {\hbox{\textbf{C}}(CR)\ni\boldsymbol{b}_c} \\
& {\hbox{\textbf{C}}(R)} \\
& {C^T\boldsymbol{\mathrm{C}}(A)=\boldsymbol{\mathrm{C}}(C^TCR) \subset \boldsymbol{\mathrm{C}}(R)}
\arrow["{{R\boldsymbol{x}_r=\boldsymbol{w}}}"{pos=0.5}, from=3-1, to=2-2]
\arrow["{CR\boldsymbol{x}_r=A\boldsymbol{x}_r=\boldsymbol{b}_c}"{pos=0.5}, from=2-2, to=3-3]
\arrow["\substack{A^TA\boldsymbol{x}_r=A^T\boldsymbol{b}_c\\[.5em]\boldsymbol{x}_r=R^T(RR^T)^{-1}R\boldsymbol{x}_r=R^+C^+\boldsymbol{b}_c}"{pos=0.3}, from=4-2, to=3-1]
\arrow[""{name=0, anchor=center, inner sep=0}, "\substack{C^TCR\boldsymbol{x}_r=C^TA\boldsymbol{x}_r=C^T\boldsymbol{b}_c\\[1ex]R\boldsymbol{x}_r=(C^TC)^{-1}C^T\boldsymbol{b}_c=C^+\boldsymbol{b}_c}"{pos=0.7}, from=3-3, to=4-2]
\arrow["{R^+}"', from=4-2, to=3-1]
\arrow["R"', from=3-1, to=2-2]
\arrow["C"', from=2-2, to=3-3]
\arrow["{C^+}"'{pos=0.45}, from=3-3, to=4-2]
\end{tikzcd}\]
\caption{Given a full rank decomposition $A=CR$, matrix $R$ maps $\hbox{\textbf{C}}(A^T)$ into $\hbox{\textbf{C}}(C^T)$ and matrix $C^T$ maps $\hbox{\textbf{C}}(A)$ into $\hbox{\textbf{C}}(R)$ when $(CR)^+=R^+C^+$. }
\label{fig:3}
\end{figure}

The necessary and sufficient conditions for $A^+ = (CR)^+ = R^+C^+$ are surprisingly complex. Greville formulated them in \cite{7}. The reverse order law for pseudoinverse, $(CR)^+=R^+C^+$, holds if and only if 
\begin{equation}
\boldsymbol{\mathrm{C}}(RR^TC^T) \subset\boldsymbol{\mathrm{C}}(C^T)
\quad\mathrm{and}\quad 
\boldsymbol{\mathrm{C}}(C^TCR) \subset \boldsymbol{\mathrm{C}}(R).
\end{equation} 
Matrix $R$ maps $\boldsymbol{x}_r = R^TC^T\boldsymbol{y}$ from the row space $\hbox{\textbf{C}}(A^T)$ into the row space $\hbox{\textbf{C}}(C^T)$. Then matrix $C$ takes $R\boldsymbol{x}_r$ to the column space $\hbox{\textbf{C}}(A)$, so we have $CR\boldsymbol{x}_r = A\boldsymbol{x}_r = \boldsymbol{b}_c$. Matrix $C^T$ maps $\boldsymbol{b}_c$ into the column space $\hbox{\textbf{C}}(R)$. Then $R^TC^T\boldsymbol{b}_c$ is also in the row space $\hbox{\textbf{C}}(A^T)$. The inverse map $C^+$ identifies $R\boldsymbol{x}_r$ in the column space $\hbox{\textbf{C}}(R)$ and then $R^+$ gives $\boldsymbol{x}_r = R^+C^+\boldsymbol{b}_c$ in Figure~ \ref{fig:3}. The reverse order law demands that 
\begin{align}
C^+C(RA^T)=RA^T 
\quad\text{and}\quad
RR^+(C^TA)=C^TA.
\end{align}
It follows that $C$ and $R$ solve the two-sided projection equation:
\begin{align}
C^+C(RR^TC^TC)RR^+ = RR^TC^TC.
\end{align}
For the full-rank factorization of Theorem~\ref{thm1} the equation holds with $C^+C = I_r = RR^+$.

\begin{figure}[!t]
\centering
\[\begin{tikzcd}
& {\hbox{\textbf{C}}((CRR^+)^+) = RR^+ \hbox{\textbf{C}}(C^T)} \\[10pt]
{[10pt]
{\boldsymbol{x}_r\in\hbox{\textbf{C}}(R^TC^T)}} & {\hbox{\textbf{C}}(R)\cap\hbox{\textbf{C}}(C^T)} & {\hbox{\textbf{C}}(CR)\ni\boldsymbol{b}_c} \\
& {\hbox{\textbf{C}}(A^T)=\hbox{\textbf{C}}((C^+CR)^+)}
\arrow["{{R\boldsymbol{x}_r = (CRR^+)^+\boldsymbol{b}_c}}"', curve={height=18pt}, from=2-3, to=2-2]
\arrow["{{\boldsymbol{x}_r=(C^+CR)^+R\boldsymbol{x}_r}}"', curve={height=18pt}, from=2-2, to=2-1]
\end{tikzcd}\]
\caption{The pseudoinverse $A^+=(CR)^+=(C^+CR)^+(CRR^+)^+$ decomposes into the product of the pseudoinverse of $R$ projected on the row space of $C^T$ and the pseudoinverse of $C^T$ projected on the column space of $R$.}
\label{fig:4}
\end{figure}
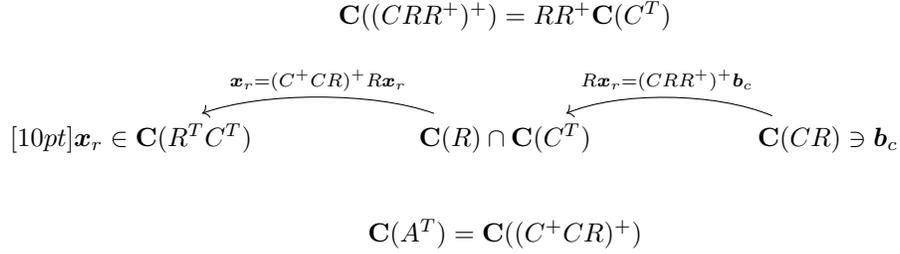

In general, the inverse map $A^+=(CR)^+$ takes $\boldsymbol{b}_c$ in the column space $\hbox{\textbf{C}}(C)$ and projects it into $R\boldsymbol{x}_r$ in the part of the column space $\hbox{\textbf{C}}(R)$ intersecting with $\hbox{\textbf{C}}(C^T)$. That is ensured by $(C^+RR^+)^+$, since by the first principles
$$
\hbox{\textbf{C}}((CRR^+)^+) = \hbox{\textbf{C}}((CRR^+)^T) = RR^+ \hbox{\textbf{C}}(C^T) = RR^+ \hbox{\textbf{C}}(C^+).
$$
Then $(C^+CR)^+$ maps $R\boldsymbol{x}_r$ into the row space  $\hbox{\textbf{C}}(R^TC^T)$. We have
$$
\hbox{\textbf{C}}((C^+CR)^+) = \hbox{\textbf{C}}((C^+CR)^T) = R^T \hbox{\textbf{C}}(C^+C) = R^T \hbox{\textbf{C}}(C^T) = \hbox{\textbf{C}}(R^TC^T).
$$
Therefore, $\hbox{\textbf{C}}(R^TC^T) = \hbox{\textbf{C}}(A^+)$ in Figure~\ref{fig:4}. Also, for any $\boldsymbol{b}_n$ in the nullspace of $R^TC^T$ we have $RR^+C^T\boldsymbol{b}_n = \boldsymbol{0}$, which shows that $\boldsymbol{b}_n$ is in the nullspace 
$$
\hbox{\textbf{N}}(RR^+C^T) = \hbox{\textbf{N}}((CRR^+)^+)\subset \hbox{\textbf{N}}((C^+CR)^+(CRR^+)^+) = \hbox{\textbf{N}}(A^+).
$$
Conversely, any $\boldsymbol{b}_n$ in the nullspace $\hbox{\textbf{N}}((CRR^+)^+) = RR^+\hbox{\textbf{N}}(C^T)$ is also in the nullspace of $\hbox{\textbf{N}}(C^T)$. It is orthogonal to the columns of $C$ spanning $\hbox{\textbf{C}}(A)$. Therefore, it is also in the nullspace $\hbox{\textbf{N}}(A^T)$. That proves the following formula for $A^+$:

\vspace{10pt}
\begin{theorem}\label{thm2}
The pseudoinverse $A^+$ of a product $A=CR$ is given by the product
\begin{equation}
A^+ = (CR)^+ = (C^+CR)^+(CRR^+)^+.
\end{equation}
\end{theorem}
\vspace{10pt}

The statement in the title of this paper is not generally true. But the statement of Theorem~2 corrects the mistake as the following example shows:
\begin{gather*}
C = \left[\begin{array}{rrr}1&0\\\end{array}\right] \quad 
R = \left[\begin{array}{rrr}1\\1\\\end{array}\right]\quad 
C^+C = \left[\begin{array}{rrr}1&0\\0&0\\\end{array}\right] \quad
RR^+ = \dfrac{1}{2}\left[\begin{array}{rrr}1&1\\1&1\\\end{array}\right]
\\[1em]
(C^+CR)^+(CRR^+)^+ = 
\left[\begin{array}{rrr}1&0\\\end{array}\right]
\left[\begin{array}{rrr}1\\1\\\end{array}\right]
= (CR)^+.
\end{gather*}
Theorem~\ref{thm3} gives an even more general statement (following from \cite{3}). 

\vspace{10pt}
\begin{theorem}\label{thm3}
The pseudoinverse $A^+$ of a product $A=CR$ is given by the product
\begin{equation}\label{eq:genROL}
A^+ = \left((P^TCR)^+P^TC\right)\left(RQ(CRQ)^+\right),
\end{equation}
for any $P$ and $Q$ satisfying
\begin{equation}
\mathrm{rank}(P^TA) = \mathrm{rank}(AQ) = \mathrm{rank}(A).
\end{equation}
\end{theorem}
\vspace{10pt}
\noindent
\textbf{Proof:}
When $\mathrm{rank}(P^TA) = \mathrm{rank}(AQ) = \mathrm{rank}(A)$ for $P$ and $Q$ of proper size, then 
$$
\hbox{\textbf{C}}((P^TA)^T) = \hbox{\textbf{C}}((A)^T) 
\quad\text{and}\quad
\hbox{\textbf{C}}(AQ) = \hbox{\textbf{C}}(A).
$$
In that case, $(P^TA)^+(P^TA) = A^+A$ and $(AQ)(AQ)^+=AA^+$ and (\ref{eq:genROL}) is true:
$$
A^+ = (P^TA)^+P^TAQ(AQ)^+.
$$
\vspace{5pt}

\begin{figure}[!ht]
\centering
\includegraphics[trim=1cm 6.5cm 1cm 6.5cm, clip, width=\textwidth]{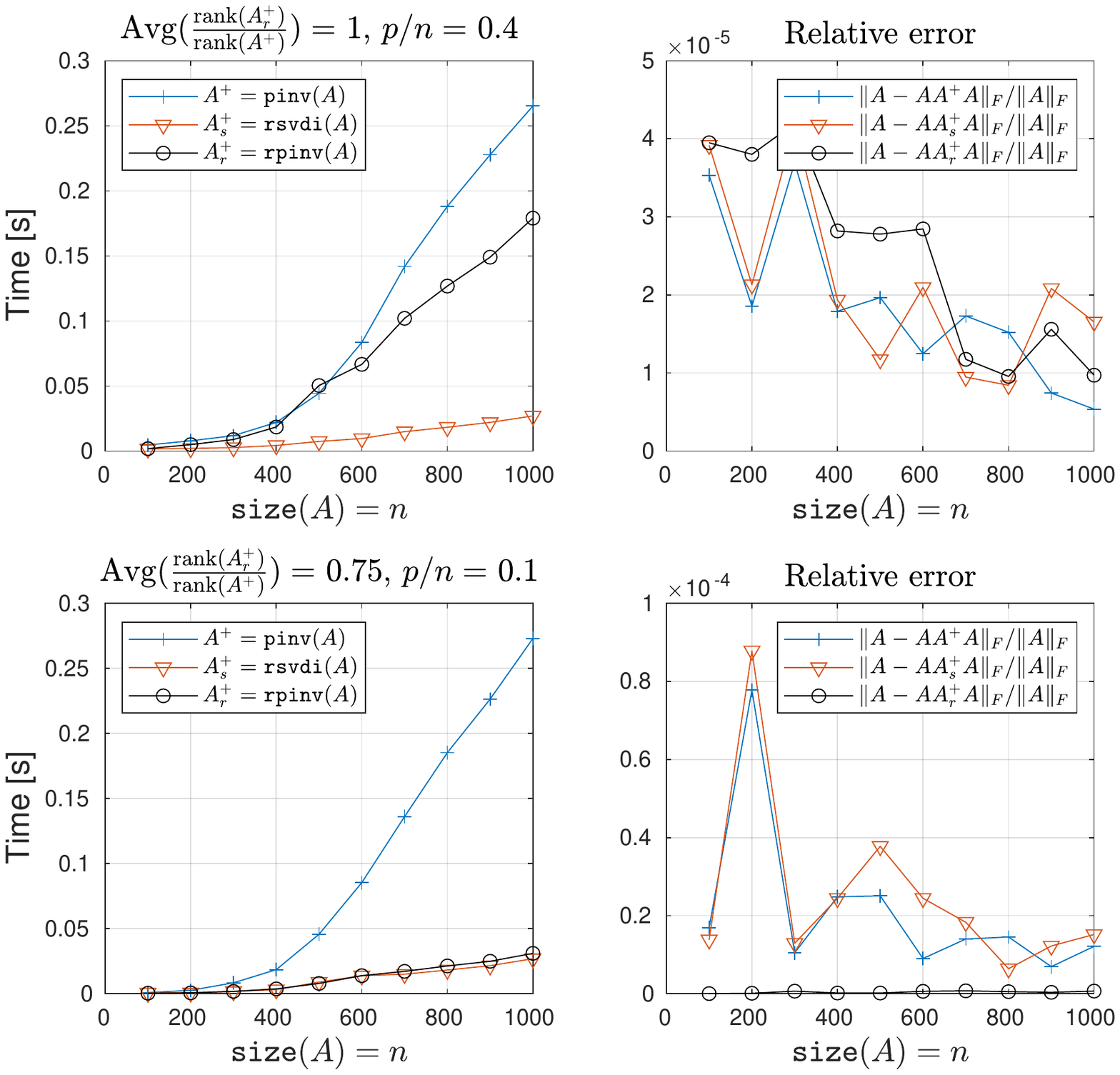}
\caption{
Computation time and relative approximation error of the randomized $A^+_r = \mathtt{rpinv}(A)$, randomized SVD-based $A^+_s = \mathtt{rsvdi}(A)$ and direct $A^+ = \mathtt{pinv}(A)$ method of calculating the pseudoinverse of $A$. In the experiment, we generate a random test matrix $A = \mathtt{gallery('randsvd',n,1e100)}$ with $n$ ranging from $100$ to $1000$. Then, we calculate the approximation $A^+_r$ taking random $m$ by $p$ matrix $P$ and $n$ by $q$ matrix $Q$ (in Theorem~\ref{thm3}) with $p = q =\lceil \alpha \cdot n\rceil$ for $\alpha = 0.4$ (top) and $\alpha = 0.1$ (bottom). Randomized SVD-based pseudoinverse is calculated based on rank $s = \mathrm{rank}(A)$ approximation $A_s = QU_s\Sigma_sV_s^T$ for $[Q,\sim]=\mathrm{qr}(A*\mathtt{randn}(n,s))$.}
\label{fig:5}
\end{figure}

Matrices $P$ and $Q$ must be rank-preserving to reconstruct $A^+$. But apart from that requirement, $P$ and $Q$ can be random matrices. The pseudoinverse of
$$
A=\left[\begin{array}{rrr}1&4&5\\2&3&5\\\end{array}\right]
\ \text{is equal to}
\
A^+=\dfrac{1}{15}\left[\begin{array}{rrr}-8&9\\7&-6\\-1&3\\\end{array}\right].
$$
Randomly selected matrices
$$
P=\left[\begin{array}{rrr}2&2&2\\1&2&2\\\end{array}\right]
\quad\text{and}\quad
Q=\left[\begin{array}{rrr}1&1\\0&2\\0&0\\\end{array}\right]
$$
preserve $\mathrm{rank}(A) = \mathrm{rank}(P^TA) = \mathrm{rank}(AQ) = 2$. Therefore, 
$$
(P^TCR)^{+}P^TC 
=\dfrac{1}{3}\left[\begin{array}{rrr}2&-1\\-1&2\\1&1\\\end{array}\right]
\quad\text{and}\quad
RQ(CRQ)^+
=\dfrac{1}{5}\left[\begin{array}{rrr}-3&4\\2&-1\\\end{array}\right]
$$
reconstruct
$$
A^+=
\dfrac{1}{3}\left[\begin{array}{rrr}2&-1\\-1&2\\1&1\\\end{array}\right]
\cdot
\dfrac{1}{5}\left[\begin{array}{rrr}-3&4\\2&-1\\\end{array}\right]
=
\dfrac{1}{15}\left[\begin{array}{rrr}-8&9\\7&-6\\-1&3\\\end{array}\right].
$$

\vspace{10pt}
Random sampling matrices $P$ and $Q$ that are not rank-preserving produce a~low-rank approximation of $A^+$ from the samples of bases for the column space and row space of $A$. Therefore, the formula in Theorem~\ref{thm3} describes a randomized algorithm approximating $A^+$:
\begin{verbatim}
function Arplus = rpinv(A,p,q)
% rpinv: randomized pseudoinverse of A 
[m,n] = size(A); 
P   = randn(m,p);
Q   = randn(n,q);
PTA = P'*A; 
AQ  = A*Q; 
Arplus = pinv(PTA)*PTA*Q*pinv(AQ);
end
\end{verbatim}
When $P$ and $Q$ are small ($m$ by $p$ and $n$ by $q$ matrices), it is more time-efficient to compute $A^+_r$ using two pseudoinverses, $(P^TA)^+$ and $(AQ)^+$, than to calculate $A^+$ directly. Also, that approach may lead to smaller relative errors (under the circumstances to be investigated) than the randomized SVD method while maintaining similar time efficiency. Figure~\ref{fig:5} demonstrates these points by comparing each method. 

We have given three formulas for the pseudoinverse of a matrix product $A = CR$ and derived them from the four fundamental subspaces of $A$. Every matrix $A$ gives an invertible map from its row space to its column space. The pseudoinverse $A^+$ inverts the row space to the column space map when the nullspace of $A^T$ and $A^+$ match. The first formula in Theorem~\ref{thm1} is correct if $C$ has independent columns and $R$ has independent rows or when Greville's conditions hold. The second in Theorem~\ref{thm2} is always correct. The third in Theorem~\ref{thm3} is correct when the sampling matrices $P$ and $Q$ preserve rank. Otherwise, it is not correct but potentially very useful, especially for large $A$.

\end{document}